\documentclass[12pt,psamsfonts]{article} 
\usepackage{a4,fullpage,amssymb,epsf,psfrag,times}
\usepackage{graphicx}
\usepackage{amsmath}
\usepackage{amsfonts}
\usepackage{enumerate}
\usepackage{latexsym}
\usepackage{epsfig}
\usepackage{amsthm}  
\usepackage{amsmath}
\usepackage{amssymb} 
\usepackage{latexsym}
\usepackage{epsfig} 
\usepackage{amssymb,latexsym}
\usepackage[all]{xy}

\newtheorem{Theorem}{Theorem}[section]

\newtheorem{lemma}[Theorem]{Lemma}
\newtheorem{definition}[Theorem]{Definition}

\newtheorem{prop}[Theorem]{Proposition}

\newtheorem{remark}[Theorem]{Remark}
\newtheorem{example}[Theorem]{Example}

\newtheorem*{namedtheorem}{\theoremname}
\newcommand{\theoremname}{testing}
\newenvironment{named}[1]{\renewcommand{\theoremname}{#1}\begin{namedtheorem}}{\end{namedtheorem}}

\DeclareMathOperator{\Lie}{Lie}
\DeclareMathOperator{\vol}{vol}
\DeclareMathOperator{\Aut}{Aut}
\DeclareMathOperator{\GL}{GL}
\DeclareMathOperator{\AGL}{AGL}
\DeclareMathOperator{\ConvHull}{ConvHull}
\DeclareMathOperator{\length}{length}

\title{Maximal ball packings of symplectic\--toric manifolds}
\author{Alvaro Pelayo and Benjamin Schmidt}

\begin{document}
\maketitle

\begin{abstract}
Let $(M, \, \sigma, \, \psi)$ be a symplectic\--toric manifold of dimension at least four.  This paper investigates the so called \emph{symplectic ball packing problem} 
in the toral equivariant setting.  We show that 
the set of toric symplectic ball packings of $M$ admits the structure of a convex polytope.  Previous work of the first author shows that up to equivalence, only $(\mathbb{CP}^1)^2$
and $\mathbb{CP}^2$ admit density one packings when $n=2$ and only $\mathbb{CP}^n$ admits density one packings when $n>2$.  In contrast, we show that
for a fixed $n\geq2$ and each $\delta \in (0,\, 1)$, there are uncountably many inequivalent 
$2n$\--dimensional symplectic\--toric manifolds with a maximal toric packing of density $\delta$.  This result follows from a general analysis of how the densities of maximal packings change while varying a given symplectic-toric manifold through a family of symplectic\--toric manifolds that are equivariantly diffeomorphic but not equivariantly symplectomorphic. \end{abstract}

\section{Introduction and Main Result}
\vskip 5pt

\subsection{Motivation}
What portion of a manifold can be filled by disjointly embedded balls?  Answers to this ball packing question depend considerably on the manifold being packed, the types of embeddings allowed, the number of balls allowed, and the radii of balls allowed.   In this paper, we consider a packing problem for a specific class of compact symplectic manifolds.  Throughout, we let $(M^{2n},\, \sigma)$ denote a 
$2n$\--dimensional compact, connected, smooth manifold equipped with a symplectic form $\sigma$.  Before stating our main result, we discuss some of the foundational results concerning ball packings of symplectic manifolds.

In \cite{Gr}, M. Gromov considered the problem of finding the supremum of the densities $$\nu_{N}
(M,\, \sigma)=\sup_{r>0} \{ \frac{N \vol(\mathbb{B}_{r})}{\vol_{\sigma}(M)}\}$$ of packings of a symplectic manifold $(M^{2n},\, \sigma)$ by a fixed number $N$ of disjoint symplectic embeddings of the standard open radius $r$ symplectic ball $\mathbb{B}_{r} \subset \mathbb{C}^n$.   Therein, Gromov found obstructions to a full (also called perfect) packing by too few balls of the same radius; specifically, he established that $\nu_{N}(\mathbb{B}_{1})\le N/2^n$ for each $1<N<2^n$.  Motivated by techniques from algebraic geometry, D. McDuff and L. Polterovich established a correspondence between symplectically embedded balls and symplectic blowing-up in \cite{McPo}.  They showed, for instance, that when $N$ is a square, $\mathbb{CP}^2$ is fully packable by $N$ disjoint symplectically embedded balls. In contrast, they also established that $\nu_{N}(\mathbb{CP}^2)<1$ for nonsquare $1<N \le 8$.    
In \cite{Bi1}, P. Biran demonstrated that for any symplectic four manifold $(M^4,\, \sigma)$ with \textit{rational} cohomology class $[\sigma] \in \textup{H}^2(M,\, \mathbb{Q})$, $\nu_{N}(M,\, \sigma)=1$ for all sufficiently large $N$.  We refer the reader to \cite{Bi2} for more in this direction. 

In this article, we consider the symplectic packing problem in a toral equivariant setting.   Specifically, we consider \textit{toric ball packings} of \textit{symplectic\--toric manifolds} (see 
Definition \ref{toricpacking}).  In this particular ball packing problem, the Euler characteristic $\chi(M)$ gives an upper bound on the number of disjoint balls in any possible packing.  For this reason, we allow the balls in an equivariant packing to have varying radii as opposed to the packings described in the previous paragraph where all balls have the same radius.  

The study of toric ball packings of symplectic\--toric manifolds was initiated by the first author in \cite{Pe1}  where the homotopy type of the space of equivariant and symplectic embeddings of a fixed ball is described, and in \cite{Pe2} where the symplectic\--toric manifolds admitting a full toric ball packing are classified.  The present work was motivated by the latter:

\begin{Theorem}\textup{\cite[Thm. 1.7]{Pe2}} \label{Alvaro}
A symplectic\--toric manifold $(M^{2n},\, \sigma,\, \psi)$ admits a full toric ball packing if and only if there exists $\lambda>0$  such that

\begin{itemize}
\item
if $n=2$, $(M^4,\,\sigma,\, \psi)$ is equivariantly symplectomorphic to either $(\mathbb{CP}^{2}, \,\lambda \cdot \sigma_{\textup{FS}})$  or a product $(\mathbb{CP}^1  \times \mathbb{CP}^1, \, \lambda \cdot(\sigma_{\textup{FS}} \oplus \sigma_{\textup{FS}}))$ (where $\sigma_{\textup{FS}}$ denotes the Fubini\--Study form and these manifolds are equipped with the standard actions of $\mathbb{T}^2$), or
\item
if $n=1$ or $n >2$, $(M^{2n},\, \sigma,\, \psi)$ is equivariantly symplectomorphic to $(\mathbb{CP}^{n},\, \lambda \cdot \sigma_{\textup{FS}})$ (where $\sigma_{\textup{FS}}$ denotes the Fubini\--Study form and this manifold is equipped with the standard action of $\mathbb{T}^n$). 
\end{itemize}
\end{Theorem}

\vskip 5pt

\subsection{Main Results}
Our motivation for this paper was to understand when it is possible to list (up to equivalence) all symplectic\--toric manifolds  admitting a maximal toric ball packing of a specified density, such as in Theorem \ref{Alvaro}.  As it turns out, this is possible only when that density is one.  To be more precise, let $\mathcal{S}^{2n}$ denote the set of equivalence classes of $2n$\--dimensional
symplectic\--toric manifolds (see Section 2 for definitions).  The maximal density function $$\Omega_{2n}: \mathcal{S}^{2n} \rightarrow (0,\,1]$$ associates to each equivalence class $[(M^{2n}, \, 
\sigma, \, \psi)]$ the largest density from all possible symplectic\--toric packings of $M$. 

Theorem \ref{Alvaro} classifies the set $\Omega_{2n}^{-1}(\{1\})$.  
In contrast with Theorem \ref{Alvaro}, we prove the following:

\begin{Theorem}\label{contrast Theorem}
Let $\mathcal{S}^{2n}$ denote the set of equivalence classes of $2n$\--dimensional symplectic\--toric manifolds and let $\Omega_{2n}: \mathcal{S}^{2n} \rightarrow (0,\,1]$ be the maximal
density function. Then $\Omega^{-1}(\{x\})$ is uncountable for all $x \in (0,1)$.
\end{Theorem}

Theorem \ref{contrast Theorem} answers \cite[Quest.\,5.1]{Pe2}. The proof of Theorem \ref{contrast Theorem} follows from the proof of the next theorem, which asserts that there are uncountable families of equivariantly diffeomorphic symplectic toric manifolds with the same maximal density that are not equivariantly symplectomorphic.

\vskip 5pt

\begin{Theorem} \label{main Theorem}
Let $\mathcal{S}^{2n}$ denote the set of equivalence classes of $2n$\--dimensional symplectic\--toric manifolds and let 
$$
\Omega_{2n}: \mathcal{S}^{2n} \rightarrow (0,\,1]
$$ 
be the maximal density function, $n\geq2$. 
Suppose that $(M, \, \sigma,\,\psi)$ is a symplectic\--toric manifold with 
Euler characteristic $\chi(M)\geq\lfloor (n+2)/2 \rfloor \cdot \lceil (n+2)/2 \rceil +1$. 
Then for any $\epsilon>0$,  
there exists a constant $c>0$ and a family $\mathcal{F}$ of equivariantly diffeomorphic 
symplectic\--toric manifolds satisfying

\begin{itemize}

\item 
$|\vol_{\sigma'}(M')-\vol_{\sigma}(M)| < \epsilon$ for all $(M', \, \sigma', \, \psi') \in \mathcal{F}$,

\item 
$\Omega_{2n}^{-1}(\{x\}) \cap \mathcal{F}$ is uncountable for all $x \in (\delta-c,\, \delta)$ 
or for all $x \in (\delta, \, \delta+c)$. 

\end{itemize}

\end{Theorem}

\vskip 5pt

To prove Theorem \ref{main Theorem} we exploit the following:

\vskip 5pt

\begin{prop} \label{convex prop}
Let $(M^{2n},\, \sigma, \, \psi)$ be a symplectic\--toric manifold of dimension $2n \geq 4$. 
The set of symplectic\--toric packings of $(M^{2n}, \, \sigma, \, \psi)$ has the structure of
a convex polytope in $\mathbb{R}^V$, where $V$ denotes the Euler characteristic $\chi(M)$
of $M$.  Moreover, $(M^{2n},\, \sigma, \, \psi)$ admits finitely many maximal toric ball packings. 
\end{prop}

The proof of Theorem \ref{main Theorem} essentially follows from the fact that there are certain linear perturbations of a given symplectic\--toric manifold along which the maximal density function is locally convex.  The proof of its local convexity follows from the Brunn-Minkowski inequality once an explicit description of the maximal density function in terms of the polytope from Proposition \ref{convex prop} is given.

\vskip 5pt

In the next section, we collect together preliminary material and reduce Theorem \ref{main Theorem} to a proposition concerning polytopes.  In the last section we prove these propositions using basic convexity techniques. 
\vskip 5pt

\centerline {\bf Acknowledgements:}  

\vskip 5pt
We are indebted to Professor Novik for bringing to our attention an explicit example showing that the density function is surjective during the first author's visit to the University of Washington, as well
as for providing us with the statement and proof of Lemma \ref{novik}.  The first author is grateful
to her for stimulating conversations and for hospitality at the University of Washington.
Part of the work of the first author was funded by Rackham Fellowships from the University of Michigan. His research was partially conducted while visiting Oberlin College.  The second author was partially funded by the Clay Mathematics Institute as a Liftoff Fellow as well as by an NSF Postdoctoral Fellowship during the period this research was conducted.

\section{Reduction to Convex Geometry}

\vskip 5pt

Throughout this section $(M^{2n}, \,\sigma)$ denotes a $2n$\--dimensional compact connected smooth manifold equipped with a symplectic form.  We let $\mathbb{T}^k\cong (S^1)^k$ denote the $k$-dimensional torus and identify its Lie algebra $\frak{t}$ with $\mathbb{R}^k$.  This identification is not unique and throughout we make the convention that the identification comes from the isomorphism $\Lie(S^1)\cong \mathbb{R}$ given by $\frac{\partial}{\partial \theta} \mapsto 1/2$.    Using the standard inner product on $\mathbb{R}^k$, we identify the dual space $\frak{t}^*$ with $\frak{t}$.

\begin{definition} \label{ham}
A $\sigma$\--preserving action $\psi: \mathbb{T}^k \times M \rightarrow M$ of a 
$k$\--dimensional torus is \textup{Hamiltonian} if for each $\xi \in \mathfrak{t}$ there exists
a smooth function $\mu_{\xi} \colon M \to \mathbb{R}$ such that
$\textup{i}_{\xi_M}\sigma=\textup{d} \mu_{\xi}$, and the map $$\mathfrak{t}
\to \textup{C}^{\infty}(M, \, \mathbb{R}),\,\,\, \xi \mapsto \mu_{\xi},$$ is a Lie algebra
homomorphism.
Here, $\xi_M$ is the vector field on $M$ infinitesimally generating the one parameter action coming from $\xi$ and the Lie algebra structure on $\textup{C}^{\infty}(M, \, \mathbb{R})$ is given by the Poisson bracket.
\end{definition}

It follows from Definition \ref{ham}
that a $\sigma$\--preserving action $\psi: \mathbb{T}^k \times M \rightarrow M$ of a 
$k$\--dimensional torus is Hamiltonian if and only if there exists a \textit{momentum map} $\mu^{M}: M \rightarrow \frak{t}^*$ satisfying Hamilton's equation $$\textup{i}_{\xi_M} \sigma=\textup{d} \langle \mu^M,\, \xi \rangle,$$ for all $\xi \in \frak{t}$. Note that a momentum map is well defined up to translation by an element of $\frak{t}^*$.  Nevertheless, we will ignore this ambiguity and refer to \textit{the momentum map}.  It is well known that if $(M^{2n},\,\sigma)$ admits an effective and Hamiltonian action of $\mathbb{T}^k$, then $k\leq n$ (see for instance \cite[Theorem 27.3]{Ca}).  The maximal case, usually referred to as a \textit{symplectic\--toric manifold} or \textit{Delzant manifold}, is a triple $(M^{2n},\,\sigma, \,\psi)$ consisting of a compact connected symplectic manifold $(M,\, \sigma)$ equipped with an effective and Hamiltonian action $\psi: \mathbb{T}^n \times M \rightarrow M$.

\begin{example} \label{stsimplex}
\normalfont
Equip the open radius $r$ ball $\mathbb{B}_{r} \subset \mathbb{C}^n$ with the standard symplectic form $\sigma_0=\frac{i}{2}\, \sum_{j} \textup{d}z_j \wedge 
\textup{d}\overline{z_j}$.  The action $\textup{Rot}: \mathbb{T}^n \times \mathbb{B}_{r}\rightarrow \mathbb{B}_{r}$ of $\mathbb{T}^n$ given by $(\theta_1,\ldots, \theta_n)\cdot(z_1,\ldots, z_n)=
(\theta_1 \, z_1, \ldots ,\theta_n \, z_n)$ is Hamiltonian.  Its momentum map $\mu^{\mathbb{B}_{r}}$ has components $\mu^{\mathbb{B}_{r}}_k=|z_k|^{2}$.  Its image, which we shall denote by $\Delta^n(r)$, is given by  $$\Delta^n(r)=\ConvHull(0,r^{2}\, e_1,\ldots, r^{2}\, e_n)\setminus \ConvHull(r^{2}\, e_1,\ldots, r^{2}\, e_n),$$ where $\{e_i\}_{i=1}^n$ is the standard basis of $\mathbb{R}^n$.  When the dimension $n$ is clear from the context, we shall write
$\Delta(r)=\Delta^n(r)$.
$\oslash$
\end{example}

The next two definitions define the packings considered in this paper.  These definitions first appeared in \cite{Pe2} in a slightly different but equivalent form.  

\begin{definition}
\textup{Let $(M^{2n},\, \sigma, \, \psi)$ be a symplectic\--toric manifold,
let $\Lambda \in \Aut(\mathbb{T}^n)$ and let $r>0$. 
A subset $B\subset M$} is said to be a $\Lambda$\--equivariantly embedded symplectic ball 
of radius $r$ \textup{if there exists a symplectic embedding $f:\mathbb{B}_{r} \rightarrow M$ with image $B$ and such that the following diagram commutes}:
\begin{eqnarray} 
\xymatrix{ \ar @{} [dr] |{\circlearrowleft}
\mathbb{T}^{n} \times \mathbb{B}_{r}  \ar[r]^{\Lambda \times f}      \ar[d]^{ \textup{Rot} }  &  \mathbb{T}^{n} \times M 
                 \ar[d]^{\psi}   \\
                   \mathbb{B}_{r}  \ar[r]^f   &       M    }. \nonumber
\end{eqnarray} 
\textup{We say that} the $\Lambda$\--equivariantly embedded symplectic ball $B$ has center 
$f(0)\in M$.

\textup{We shall say that a 
subset $B'\subset M$} is an equivariantly embedded symplectic ball
of radius $r'$ \textup{if there exists $\Lambda' \in \Aut(\mathbb{T}^n)$ such that
$B'$ is a $\Lambda'$\--equivariantly embedded symplectic ball of radius $r'$ (although this
is a slight abuse of the standard use of the word ``equivariantly'')}.  
$\oslash$
\end{definition}

We define the symplectic volume of a subset $A\subset M$ by $\vol_{\sigma}(A):=\int_{A} \sigma^{n}$. 

\begin{definition} \label{toricpacking}
\textup{Let $(M^{2n},\,\sigma,\, \psi)$ be a symplectic\--toric manifold.  A} toric ball packing of $M$ \textup{is a disjoint union $\mathcal{P}:=\bigsqcup_{\alpha \in A} B_{\alpha}$ of equivariantly embedded symplectic balls $B_{\alpha}$ (of possibly varying radii) in $M$.  The} 
density $\Omega(\mathcal{P})$ of a packing $\mathcal{P}$ \textup{is defined by 
$\Omega(\mathcal{P}):= \vol_{\sigma}(\mathcal{P})/\vol_{\sigma}(M)$. The} density 
$\Omega(M^{2n},\,\sigma,\,\psi)$ of a symplectic\--toric manifold $(M^{2n},\,\sigma,\, \psi)$ \textup{is defined by $$\Omega(M,\, \sigma,\, \psi):= \sup \{\, \Omega(\mathcal{P})\, | \, \mathcal{P} \, \textup{is a toric packing of}\, M \}.$$  A packing achieving this density is said to be a} maximal density packing.  \textup{If in addition, this density is one, then $(M^{2n},\,\sigma,\, \psi)$ is said to admit a} full \textup{or} perfect toric ball packing. $\oslash$
\end{definition}

\vskip 5pt

For a symplectic\--toric manifold $(M^{2n},\, \sigma,\, \psi)$, the number of fixed points of $\psi$ is known to coincide with the Euler characteristic $\chi(M)$ (see e.g. \cite{Fu}).  It follows that a toric ball packing $\mathcal{P}$ consists of at most $\chi(M)$ disjoint equivariant balls.
By a well known theorem of Atiyah and Guillemin\--Sternberg, the image of the momentum map of a toral Hamiltonian action is the convex hull of the images of the fixed points (see for instance \cite[Theorem 27.1]{Ca}).  The images of momentum maps for symplectic\--toric manifolds are a particular class of polytopes. Recall that an $n$\--dimensional polytope is \emph{simple} if 
there are precisely $n$ edges meeting at each one of its vertices.

\begin{definition} \label{delzantdef}
\textup{A simple $n$\--dimensional convex polytope $\Delta \subset \mathbb{R}^n$ is said to be} Delzant \textup{if for each vertex $v$, the edges meeting at $v$ are all of the form $v+t_i\,u_i$ where $t_i>0$ and $\{u_1,\ldots, u_n\}$ define a basis of the $\mathbb{Z}$\--module $\mathbb{Z}^n$}. $\oslash$
\end{definition}

A polytope is describable as the intersection of closed half\--spaces $$\Delta:= \bigcap_{i=1}^{F} \{x \in \mathbb{R}^n \, \vert \,  \langle x,\, u_i \rangle \geq \lambda_i \},$$ where the vector $u_i$ is an inward pointing normal vector to the $i^{\textup{th}}$ facet of $\Delta$ and each $\lambda_i$ is a real scalar.  In this notation, the polytope $\Delta$ is Delzant if and only if there are precisely $n$ facets incident to each vertex of $\Delta$ and the inward pointing normals to these facets $u_1, \ldots, u_n$ can be chosen to be a $\mathbb{Z}$\--basis of $\mathbb{Z}^n$.  

\vskip 5pt
Given a symplectic\--toric manifold $M$, the symplectic\--toric manifolds obtained from $M$ by scaling the symplectic form, changing the time parameter in the acting torus by an automorphism, and any others that are equivariantly symplectomorphic to one of these will clearly have the same maximal density.  Therefore, for the purpose of this paper we will say 
that \textit{two symplectic\--toric manifolds $(M_1,\, \sigma_1,\, \psi_1)$ and 
$(M_2,\, \sigma_2,\, \psi_2)$ 
are equivalent} if there exists an 
automorphism $\Lambda \in \Aut(\mathbb{T}^n)\cong \GL(n,\, \mathbb{Z})$, a 
positive number $\lambda>0$, and a 
symplectomorphism $h:(M_1,\, \sigma_1) \rightarrow (M_2, \, \lambda \cdot \sigma_2)$ such that the following diagram commutes:
\begin{eqnarray} 
\xymatrix{ \ar @{} [dr] |{\circlearrowleft}
\mathbb{T}^{n} \times M_1  \ar[r]^{\Lambda \times h}      \ar[d]^{\psi_1}  &  \mathbb{T}^{n} \times M_2 
                 \ar[d]^{\psi_2}   \\
                   M_1  \ar[r]^h   &       M_2    }. \nonumber
\end{eqnarray} 

We recall the result of Delzant \cite{De}:

\begin{Theorem}\label{Delzant}
Suppose that $(M_1,\, \sigma_1,\, \psi_1)$ and $(M_2,\, \sigma_2,\, \psi_2)$ are two $2n$\--dimensional symplectic\--toric manifolds with momentum maps $\mu^{M_1}$ and $\mu^{M_2}$ respectively.  Then there exists a $(\psi_1,\,\psi_2)$\--equivariant symplectomorphism $h:(M_1,\,\sigma_1) \rightarrow (M_2,\,\sigma_2)$  such that $\mu^{M_1}=\mu^{M_2} \circ h$ if and only if 
$\mu^{M_1}(M_1)=\mu^{M_2}(M_2).$
\end{Theorem}

In view of this theorem, there is a natural equivalence relation one can put on the set of Delzant polytopes so that momentum maps will induce a bijective correspondence between equivalence classes of symplectic\--toric manifolds as above and equivalence classes of Delzant polytopes.  To be more specific, first note that scaling $\mathbb{R}^n$ leaves invariant the Delzant 
polytopes.  Define two Delzant polytopes $\Delta_1,\Delta_2\subset \mathbb{R}^n$ \textit{to be in the 
same projective class} if there exists $\lambda>0$ such that $\lambda \, \Delta_1=\Delta_2$.  The group $\AGL(n,\, \mathbb{Z})$ consisting of affine 
transformations of $\mathbb{R}^n$ with linear part in $\GL(n,\, \mathbb{Z})$ acts 
on the set of  projective classes of Delzant polytopes in $\mathbb{R}^n$.  For the
purpose of this paper we say that \textit{two Delzant polytopes $\Delta_1$ and $\Delta_2$ are equivalent} if the projective classes of $\Delta_1$ and $\Delta_2$ are in the same $\AGL(n,\, \mathbb{Z})$ orbit.  By applying Theorem \ref{Delzant}, it is standard to show that there is a bijective correspondence between equivalence classes of symplectic\--toric manifolds and equivalence classes of Delzant polytopes as defined here.

We exploit this correspondence in order to reduce Theorem \ref{main Theorem} and Proposition \ref{convex prop} to propositions concerning packing Delzant polytopes.  The next few definitions are translations of the above definitions into the appropriate definitions concerning Delzant polytopes.  Again, they first appeared in \cite{Pe2} in a slightly different but equivalent form.

\begin{definition} \label{adsi}
\textup{Let $\Delta$ be a Delzant polytope.  A subset $\Sigma \subset \Delta$ is said to be an} admissible simplex of radius $r$ with center at the vertex $v\in \Delta$ \textup{if $\Sigma$ is the image of $\Delta(r^{1/2})$ by an element of $\AGL(n,\,\mathbb{Z})$ which takes the origin to $v$ and the edges of $\Delta(r^{1/2})$ to the edges of $\Delta$ meeting $v$.  For a vertex $v \in \Delta$, we put} 
\begin{eqnarray}
r_v:=\max\{r>0\, |\, \exists \, \textup{an admissible simplex of radius $r$ with center $v$} \}.\hspace{7mm}\oslash
\label{rv}
\end{eqnarray}

\end{definition}

\begin{remark}
\normalfont
In view of Example \ref{stsimplex}, the simplex $\Delta(r^{1/2})$ may be identified with the set obtained by removing from $\ConvHull(0,re_1, \ldots, re_n)$ the facet not containing the origin.  For this reason we say that $\AGL(n,\,\mathbb{Z})$ images of $\Delta(r^{1/2})$ as in the above definition have radius $r$ instead of radius $r^{1/2}$.  $\oslash$.
\end{remark}

We denote the Euclidean volume of a subset $A\subset \Delta$ by $\vol_{\textup{euc}}(A)$.

\begin{definition}
\textup{Let $\Delta$ be a Delzant polytope.  An} admissible packing of $\Delta$ \textup{is a disjoint union $\mathcal{P}:=\bigsqcup_{\alpha \in A} \Sigma_{\alpha}$ of admissible simplices (of possibly varying radii) in $\Delta$.  The} density $\Omega(\mathcal{P})$ of a packing $\mathcal{P}$ \textup{is defined by  $\Omega(\mathcal{P}):= \vol_{\textup{euc}}(\mathcal{P})/\vol_{\textup{euc}}(\Delta)$.  The} density $\Omega(\Delta)$ of a Delzant polytope $\Delta$ \textup{is defined by $$\Omega(\Delta):= \sup \{\, \Omega(\mathcal{P})\, | \, \mathcal{P} \, \textup{is an admissible packing of}\, \Delta \}.$$  A packing achieving this density is said to be a} maximal density packing.  \textup{If in addition, this density is one, then $\Delta$ is said to admit a} full or perfect packing. $\oslash$
 \end{definition}
 
 \vskip 5pt

The next lemma shows that admissible simplices in $\Delta$ are parametrized by their centers and radii.  The \textit{rational or $\textup{SL}(n,\,\mathbb{Z})$\--length} of an interval $I\subset \mathbb{R}^n$ with rational slope is the unique number $l:=\length_{\mathbb{Q}}(I)$ such that $I$ is $\AGL(n,\, \mathbb{Z})$-congruent to an interval of length $l$ on a coordinate axis.  For a vertex $v$ in a Delzant polytope, we denote the $n$ edges leaving $v$ by $e_v^1,\ldots, e_v^n$.  By Definition \ref{delzantdef}, each $e_v^i$ is of the form $v+t_v^i \,u_v^i$ with $t_v^i>0$ and $\{u_v^i\}_{i=1}^{n}$ defining a $\mathbb{Z}$\--basis of $\mathbb{Z}^n$.  In this notation, we have that $\length_{\mathbb{Q}}(e_v^i)=t_v^i$.

 \begin{lemma} \label{adsi}
 Let $\Delta$ be a Delzant polytope.  Then for each vertex $v\in \Delta$, 
 $$r_v=\min\{\length_{\mathbb{Q}}(e_{v}^1), \ldots, \length_{\mathbb{Q}}(e_{v}^n) \}.$$  
 There is an admissible simplex $\Sigma(v,\,r)$ of radius $r$ with center $v$
 if and only if $0\leq r \leq r_v$.  Moreover this admissible simplex is unique
 and $\vol_{\textup{euc}}(\Sigma(v,\, r))= r^{n}/n!.$
 \end{lemma}
 
\begin{proof}
We first argue the uniqueness of an admissible simplex $\Sigma$ of radius $r$ and center $v$, assuming its existence.  Suppose $\Sigma_1$ and $\Sigma_2$ were two such.  By definition, there exist two affine transformations $A_1,\, A_2 \in \AGL(n,\,\mathbb{Z})$ satisfying $A_1(\Delta(r^{1/2}))=\Sigma_1$ and $A_2(\Delta(r^{1/2}))=\Sigma_2$.  As both $\Sigma_i$ are centered at $v$, both $A_i$ have translational part given by $v$.  Write $A_i(\cdot)=\Lambda_i(\cdot) +v$ with $\Lambda_i \in 
\GL(n,\, \mathbb{Z})$.  By the Delzant property of $\Delta$, c.f. Definition \ref{delzantdef}, the automorphisms $\Lambda_1$ and $\Lambda_2$ both take the standard basis $\{e_i\}_{i=1}^{n}$ of $\mathbb{Z}^n$ bijectively onto the basis $\{u_v^i\}_{i=1}^{n}$ of $\mathbb{Z}^n$ as unordered sets.  Therefore $\Lambda_{1}^{-1} \Lambda_2$ leaves invariant the standard basis of $\mathbb{Z}^n$ as an unordered set and hence leaves $\Delta(r^{1/2})$ invariant as a set.  It follows that $\Sigma_1=\Sigma_2$. Next we argue existence.

As above, write $e_{v}^i=v+t_v^i \, u_v^i$ with $\{u_v^i\}_{i=1}^{n}$ forming a $\mathbb{Z}$\--basis of $\mathbb{Z}^n$.  Let $\Lambda \in \GL(n,\,\mathbb{Z})$ be the automorphism of $\mathbb{Z}^n$ defined by $\Lambda(e_i)=u_v^i$, with $\{e_i\}_{i=1}^{n}$ the standard basis of $\mathbb{Z}^n$.  The affine transformation $A:\mathbb{R}^n \rightarrow \mathbb{R}^n$ defined by $A(x)=\Lambda(x) + v$ satisfies $A(t\, e_{i})=t\, u_v^i+v$ for each $t>0$.  Therefore, $A(\Delta(r^{1/2}))$ defines an admissible simplex for $\Delta$ with center $v$ if and only if $r\leq \min\{t_v^1, \ldots, t_v^n \}=\min\{\length_{\mathbb{Q}}(e_{v}^1), \ldots, \length_{\mathbb{Q}}(e_{v}^n) \}$, concluding the proof of existence.

Finally, for each $0\leq r \leq r_v$, let $\Sigma(v,\, r)$ be the unique admissible simplex with radius $r$ and center $v$ guaranteed by the previous two paragraphs.  
Since elements in $\AGL(n,\,\mathbb{Z})$ preserve Euclidean volume, we have 
$\vol_{\textup{euc}}(\Sigma(v,\,r))=\vol_{\textup{euc}}(\Delta(r^{1/2}))=r^{n}/n!$, completing the proof.   
\end{proof}

\vskip 5pt

Lemmas \ref{2.11} and \ref{2.12} below reduce the problem of analyzing the densities of toric ball packings of symplectic\--toric manifolds to analyzing the densities of admissible packings of Delzant polytopes.  These lemmas appear in \cite{Pe1} and \cite{Pe2} respectively, in slightly different form.  We include the arguments here for the reader's convenience.  We first recall some definitions and a result due to Y. Karshon and S. Tolman.  
 
 Let $(N^{2n}, \, \sigma^N)$ denote a connected (and not necessarily compact) symplectic manifold equipped with an effective Hamiltonian action of $\mathbb{T}^n$ with momentum map $\mu^N:N \rightarrow \frak{t}^*$.  Suppose that $\mathcal{T}\subset \frak{t}^*$ is an open and convex subset containing $\mu^N(N)$ with the property that $\mu^N:N \rightarrow \mathcal{T}$ is a proper map.  The quadruple $(N,\,\sigma^N,\,\mu^N,\,\mathcal{T})$ is said to be a \textit{proper Hamiltonian $\mathbb{T}^{n}$-manifold}.  For a subgroup $K\subset \mathbb{T}^n$, denote the fixed point set for the $\mathbb{T}^n$ action on $N$ by $N^{K}$.  A proper Hamiltonian $\mathbb{T}^n$\--manifold $(N,\, \sigma^N,\, \mu^N, \,\mathcal{T})$ is said to be \textit{centered} at $\alpha \in  \mathcal{T}$ provided that $\alpha$ is in the momentum map image of every component of $N^{K}$ for all subgroups $K\subset \mathbb{T}^n$.  For a fixed point $p \in N^{\mathbb{T}^n}$, there exist \textit{isotropy weights} $\eta_1, \ldots, \eta_n \in \frak{t}^*$ such that the induced linear symplectic action on 
$\textup{T}_p(N)$ is isomorphic to 
 the action on $(\mathbb{C}^n, \, \sigma_0)$ generated by the momentum map $$\mu^{\mathbb{C}^n}(x)=\mu^{N}(p)+ \sum_{j} |z_j|^2 \eta_j.$$  For a symplectic\--toric manifold $M$ with fixed point $p\in M^{\mathbb{T}^n}$, and corresponding vertex $v:=\mu(p) \in \Delta$, and edges $e_v^i=v+t_v^i \, u_v^i$ ($i=1, \ldots, n$) emanating from $v$ , the $\mathbb{Z}$\--basis  $\{u_v^i\}_{i=1}^{n}$ of $\mathbb{Z}^n$ coincides with the set of isotropy weights at $p$. 
 
\begin{prop}\textup{\cite[Prop. 2.8]{KaTo}} \label{2.9}
Let $(N^{2n}, \, \sigma^N,\, \mu^N,\, \mathcal{T})$ be a proper Hamiltonian $\mathbb{T}^n$-manifold.  Assume that $N$ is centered about $\alpha \in \mathcal{T}$ and that $(\mu^N)^{-1}(\{\alpha\})$ consists of a single fixed point $p$.  Then $N$ is equivariantly symplectormorphic to $$\{z \in \mathbb{C}^n \, | \, \alpha + \sum_{j} |z_j|^{2}\eta_j \in \mathcal{T}\},$$ where $\eta_1,\ldots, \eta_n \in \frak{t}^{*}$ are the isotropy weights at $p$.
\end{prop}
 
\begin{remark}
\normalfont
The reader who consults \cite{KaTo} will find an additional factor of $\pi$ in the definition of isotropy weights and in the statement of their Proposition 2.8.  This factor does not appear in the present work because of the particular identification chosen between $\frak{t}$ and $\mathbb{R}^n$. 
$\oslash$
\end{remark}

 \begin{lemma} \label{2.11}
Let $(M^{2n}, \, \sigma,\, \psi)$ be a symplectic\--toric 
manifold with momentum map $\mu^{M}: M \rightarrow \frak{t}^{*}$ and 
associated Delzant polytope $\Delta:=\mu^M(M)$.  Let $B\subset M$ be an 
equivariantly embedded symplectic ball of radius $r$ and center $p\in M$.  Then $\mu^M(B)$ is an admissible simplex of radius $r^{2}$ in $\Delta$ with center $\mu^{M}(p)$.  Conversely, if $\Sigma \subset \Delta$ is an admissible simplex of radius $r$, then there exists an equivariantly embedded symplectic ball $B\subset M$ of radius $r^{1/2}$ satisfying $\mu^M(B)=\Sigma$.
 \end{lemma}
 
 \begin{proof}
 First suppose that $B\subset M$ is an equivariantly embedded ball of radius $r$.  By definition, there is an automorphism $\Lambda \in \Aut(\mathbb{T}^n)\cong \GL(n,\, \mathbb{Z})$ and a symplectic embedding $f:(\mathbb{B}_{r},\, \sigma_{0}) \rightarrow (M,\, \sigma)$ with image $B$ such that the following diagram commutes: 
 \begin{eqnarray} 
\xymatrix{ \ar @{} [dr] |{\circlearrowleft}
\mathbb{T}^{n} \times \mathbb{B}_{r}  \ar[r]^{\Lambda \times f}      \ar[d]^{ \textup{Rot} }  &  \mathbb{T}^{n} \times M 
                 \ar[d]^{\psi}   \\
                   \mathbb{B}_{r}  \ar[r]^f   &       M    }. \nonumber
\end{eqnarray} 
We denote by $\xi_{\mathbb{B}_{r}}$ and $\xi_{M}$ the vector fields on $\mathbb{B}_{r}$ and $M$ infinitesimally generating the actions of the one paramater group coming from $\xi \in \frak{t}$.  Fix $z \in \mathbb{B}_{r}$ and a tangent vector $v\in \textup{T}_{z} \mathbb{B}_{r}$.

Note that from the definition of the momentum maps $\mu^{M}$ and $\mu^{\mathbb{B}_{r}}$, commutativity of the above diagram, and the fact that $f^*\sigma=\sigma_0$, 
we have the following sequence of equalities: 
\begin{eqnarray}
\langle \textup{d}\mu^{\mathbb{B}_{r}}_{z} (v),\, \xi \rangle_{\mu^{\mathbb{B}_{r}}(z)} 
 &=& \sigma_{0}(v, \, \xi_{\mathbb{B}_{r}})_{z} \nonumber \\
&=&  \sigma(\textup{d}f_{z}({v}),\, \textup{d}f_{z}(\xi_{\mathbb{B}_{r}}))_{f(z)} 
\nonumber \\
&=&  \sigma(\textup{d}f_{z}({v}),(\Lambda\, \xi)_{M})_{f(z)} 
\nonumber \\
&=& \langle \textup{d}\mu^{M}_{f(z)} (\textup{d}f_{z}({v})),\, \Lambda\,\xi 
\rangle_{\mu^{M}(f(z))} \nonumber \\
&=& \langle \Lambda^{\textup{t}}\,\textup{d}\mu^{M}_{f(z)}(\textup{d}f_{z}(v)),\, \xi \rangle_{\mu^{\mathbb{B}_r}(z)}. 
 \label{eq1}
\end{eqnarray}
By equation (\ref{eq1}) and the chain rule we 
obtain that for all $ {z} \in \mathbb{B}_{r}$
and $ {v} \in \textup{T}_{z} \mathbb{B}_{r},$
\begin{eqnarray}
\textup{d}\mu^{\mathbb{B}_{r}}_{z} ( {v})=\textup{d}(\Lambda^{\textup{t}}
\circ \mu^{M} \circ f)_{z} ( {v}). \label{eq2}
\end{eqnarray}

As $\mathbb{B}_r$ is connected, (\ref{eq2}) implies that there exists $x' \in \mathbb{R}^n$ such that $\mu^{\mathbb{B}_{r}}+x' =  \Lambda^{\textup{t}}\circ \mu^{M} \circ f$ as maps $\mathbb{B}_r \rightarrow \mathbb{R}^n.$  Letting $x=(\Lambda^{\textup{t}})^{-1}(x')$ yields commutativity of the following diagram:
\begin{eqnarray}  
\xymatrix{ \ar @{} [dr] |{\circlearrowleft}
\Delta(r) \ar[r]^{(\Lambda^{\textup{t}})^{-1}+ {x}}      &  \Delta^M  \\
           \mathbb{B}_{r}  \ar[r]^f  \ar[u]^{\mu^{\mathbb{B}_{r}}}  &       M      \ar[u]_{\mu^M}}.
\end{eqnarray}
It follows from commutativity of this diagram that $x=\mu^M(f(0))=\mu^M(p)$ so that $\mu(B)=\mu(f(\mathbb{B}_{r}))$ is an admissible simplex of radius $r^{2}$ and center $\mu^{M}(p)$, completing the proof of the first statement.  

Next suppose that $\Sigma \subset \Delta$ is an admissible simplex of radius $r$.  By applying a translation, we may assume that $\Sigma$ is centered at the origin.  Identify $\Sigma$ with the set 
$$
\ConvHull(0, \,r\,\eta_1,\ldots, r\,\eta_n)\setminus\ConvHull(r\,\eta_1, \ldots, r\,\eta_n).  
$$
Let $\mathcal{T}\subset \frak{t}^{*}$ be the the unique open half space of $\frak{t}^{*}$ containing $\Sigma$ with bounding hyperplane containing  $\ConvHull(r\,\eta_1, \ldots r\,\eta_n)$.  Denote by $\sigma^{N}, \, \mu^N$ the restrictions of the symplectic form $\sigma$, and of the
momentum map $\mu^M$, to the open 
submanifold $N:=(\mu^M)^{-1}(\Sigma)\subset M$.  
The quadruple $(N,\,\sigma^{N},\,\mu^N, \,\mathcal{T})$ is a proper Hamiltonian $\mathbb{T}^n$ manifold centered at $0\in \mathcal{T}$.  It now follows from Proposition \ref{2.9}, that $(N,\, \sigma^{N})$ is equivariantly symplectomorphic to 
$$
\{z \in \mathbb{C}^n \, | \, \sum_{j} |z_j|^{2}\eta_j \in \mathcal{T}\}=\mathbb{B}_{r^{1/2}}.
$$  
In other words, the set 
$N\subset M$ is an equivariantly embedded symplectic ball of radius $r^{1/2}$ satisfying $\mu^M(N)=\Sigma$ (c.f. \cite[Lem 2.3]{Pe2} for an explicit verification).
\end{proof}

\begin{lemma} \label{2.12}
 Let $(M^{2n},\, \sigma,\, \psi)$ be a symplectic\--toric manifold with momentum map $\mu^M: M \rightarrow \mathbb{R}^n$ and associated Delzant polytope $\Delta:=\mu^M(M)$.  Then for each toric ball packing $\mathcal{P}$ of $M$, $\mu^M(\mathcal{P})$ is an admissible packing of $\Delta$ satisfying $\Omega(\mathcal{P})=\Omega(\mu^M(\mathcal{P}))$.  Moreover, given an admissible packing $\mathcal{Q}$ of $\Delta$, there exists a toric ball packing $\mathcal{P}$ of $M$ satisfying $\mu^M(\mathcal{P})=\mathcal{Q}$. 
 \end{lemma}

\begin{proof}
Let $\mathcal{P}$ be a toric ball packing of $M$.  For each equivariant symplectic ball $B$ in the packing $\mathcal{P}$, $\mu^M(B)$ is an admissible simplex in $\Delta$ by the previous lemma.  Since the fibers of the momentum map $\mu^M: M \rightarrow \Delta$ are connected (see for instance \cite[Theorem 27.1]{Ca}), disjoint equivariant symplectic balls in the packing $\mathcal{P}$ are sent to disjoint admissible simplices in $\Delta$.  Hence, $\mu^M(\mathcal{P})$ is a toric packing of $\Delta$.  By the Duistermaat\--Heckman Theorem (see for instance \cite[Theorem 30.3]{Ca}) the push forward of symplectic volume under the momentum map satisfies $\mu^M_{*}(\vol_M)=K(n) \vol_{\textup{euc}}$, where $K(n)>0$ is a dimensional constant.  It follows that $\Omega(\mathcal{P})=\Omega(\mu^M(\mathcal{P}))$.  It remains to argue that given an admissible packing $\mathcal{Q}$ of $\Delta$, there exists a toric ball packing $\mathcal{P}$ of $M$ satisfying $\mu^M(\mathcal{P})=\mathcal{Q}$.  By 
the Lemma \ref{2.11}, for each admissible simplex $\Sigma$ there exists an equivariant symplectic ball $B\subset M$ with $\mu^M(B)=\Sigma$.  Choosing one such equivariant symplectic ball for each admissible simplex in $\mathcal{Q}$ defines a disjoint collection $\mathcal{P}$ of equivariant symplectic balls mapping onto $\mathcal{Q}$ under the momentum map.  Hence $\mathcal{P}$ is a toric packing of $M$ satisfying $\mu^M(\mathcal{P})=\mathcal{Q}$.
\end{proof}

\vskip 5pt

Before concluding this section by reformulating our main results in terms of polytopes, we must first introduce the complete regular $n$\---dimensional fan associated to an $n$\--dimensional Delzant polytope 
$\Delta$, as in \cite[Sec. 5]{De}.  As above, write 
$$\Delta:= \bigcap_{i=1}^{F} \{x \in \mathbb{R}^n \, \vert \,  \langle x,\, u_i \rangle \geq \lambda_i \},$$ where $F$ is the number of facets of $\Delta$ and $u_i$ is the unique primitive integral normal vector to the $i^{\textup{th}}$ facet.  For each face $\Delta'$ of $\Delta$ of codimension $k$, there is a unique multi\--index $I_{\Delta'}$ of length $k$, $I_{\Delta'}=\{i_1,\ldots,i_k\}$, $1\leq i_1<\ldots<i_k\leq F$, such that 
$$\Delta'=\{x \in \mathbb{R}^n \, \vert \, \langle x,\, u_j \rangle= \lambda_j, \,\,\forall j \in I_{\Delta'}\}.
$$   
Letting $\sigma_{\Delta'}$ denote the cone in $\mathbb{R}^n$ generated by the vectors $\{u_j \, \vert \, j \in I_{\Delta'} \}$, the \emph{complete regular $n$\--dimensional fan associated to $\Delta$} is given by $\{ \sigma_{\Delta'} \, \vert \, \Delta' \, \textup{is a face of}\, \Delta \}$.  For our purposes here, we state the well known fact that if two Delzant polytopes have the same associated fan, their associated symplectic\--toric manifolds are equivariantly diffeomorphic. This is a standard fact that follows from the 
construction of a symplectic-toric manifold starting from a given Delzant polytope.

It now follows from Lemma \ref{2.11} and Lemma \ref{2.12} that proving Theorem \ref{contrast Theorem}  is reduced to establishing the following proposition:

\begin{prop}\label{contrast prop}
Let $\mathcal{D}^n$ denote the set of equivalence classes of $n$\--dimensional Delzant polytopes and $$\Omega:\mathcal{D}^n\rightarrow (0,\,1]$$ be the maximal density function, $n\geq 2$.  Then $\Omega^{-1}(\{x\})$ is uncountable for all $x \in (0,\,1)$. 
\end{prop}

Recall that the number of fixed points 
of the $\mathbb{T}^n$\--action on $M$ equals the Euler Characteristic of $M$
and  the number of vertices of the momentum polytope $\mu(M)$,  c.f. \cite{Fu}.
We are grateful to Professor Novik for the following observation:

\begin{lemma} \label{novik}
An $n$\--dimensional 
Delzant polytope $\Delta$ with at least
$
\lfloor (n+2)/2 \rfloor \cdot \lceil (n+2)/2 \rceil +1
$
vertices has at least $n+3$ facets. Moreover, this bound is sharp in the sense
that there exists an $n$\--dimensional Delzant polytope with
$
\lfloor (n+2)/2 \rfloor \cdot \lceil (n+2)/2 \rceil
$
vertices and $n+2$ facets.
\end{lemma}
\begin{proof}
Indeed, to see that a fewer number of vertices is not enough, let $\Delta$ be
the direct product of two 
regular simplexes, one of dimension $\lfloor n/2 \rfloor$ and another one of dimension 
$\lceil n/2 \rceil$. Their product is an $n$\--dimensional Delzant polytope that has 
$\lfloor (n+2)/2 \rfloor \cdot \lceil (n+2)/2 \rceil$ vertices and only $(n+2)$ facets.
It is well known, see e.g. \cite[pp. 98-100]{grunbaum} for the proof of (the dual)
statement, that
$\Delta$ has the maximal number of vertices amongst 
all simple $n$\--polytopes with $n+2$ facets. 
\end{proof}

It then follows from Lemma \ref{2.11}, Lemma \ref{2.12} and Lemma \ref{novik}
that proving Theorem \ref{main Theorem}
is reduced to establishing the following proposition:

\begin{prop} \label{main prop}
Let $\mathcal{D}^n$ denote the set of equivalence classes of $n$\--dimensional Delzant polytopes and $$\Omega:\mathcal{D}^n\rightarrow (0,\,1]$$ be the maximal density function, $n\geq2$.  
Suppose that $\Delta$ is a Delzant polytope having at least $n+3$ facets and let $\Omega(\Delta):=\delta \in (0,\, 1)$.
Then for any $\epsilon>0$,  
there exists a constant $c>0$ and a family $\mathcal{F}$ of Delzant polytopes satisfying 

\begin{itemize}

\item the polytopes in $\mathcal{F}$ determine a common fan, 

\item 
$|\vol_{\textup{euc}}(\Delta')-\vol_{\textup{euc}}(\Delta)| < \epsilon$ for all $\Delta'  \in \mathcal{F}$,

\item 
$\Omega^{-1}(\{x\}) \cap \mathcal{F}$ is uncountable for all $x \in (\delta-c,\, \delta)$ or for all
$x \in (\delta, \, \delta+c)$. 

\end{itemize}

\end{prop}

Similarly we have reduced proving Proposition \ref{convex prop} to showing:

\begin{prop} \label{convex prop2}
Let $\Delta \subset \mathbb{R}^n$ be a Delzant polytope. The set of admissible packings
of $\Delta$ has the structure of a convex polytope in $\mathbb{R}^V$, where $V$
is the number of vertices of $\Delta$.  Moreover, $\Delta$ admits finitely many maximal density packings.
\end{prop}

We prove these propositions in the next section.

\section{Proofs of propositions \ref{contrast prop}, \ref{main prop}, and \ref{convex prop2}.}
\vskip 5pt

In this section, we prove Propositions \ref{contrast prop}, \ref{main prop}, and \ref{convex prop2} using convexity arguments. First we recall some preliminary notions.  For a set $A\subset \mathbb{R}^n$, denote the closure of $A$ by $\overline{A}$.  Now suppose that $A$ is a convex set.  A function $f:A\rightarrow \mathbb{R}$ is said to be \emph{convex} if for all 
$x_1,\, x_2 \in A$ and $t \in (0,\, 1)$, $$f((1-t)\, x_1+t\, x_2)\leq (1-t)\, f(x_1)+t\, f(x_2)$$ 
and \emph{strictly convex} if the inequality is always strict.  Similarly, the function $f$ is said 
to be \emph{concave} if for all $x_1,\, x_2 \in A$ and $t \in (0,\,1)$, 
$$f((1-t)\, x_1+t\, x_2)\geq (1-t)\, f(x_1)+t\, f(x_2)$$ and \emph{strictly concave} if the inequality is always strict.  It follows from the definitions that if $f$ is a convex function on $A$ and $g$ is a positive concave function on $A$, then $f/g$ is a convex function which is strict if one of $f$ or $g$ is strict.  Moreover, if $f_1,\ldots,f_k$ are convex functions on $A$ then so is the function $\max\{f_1,\ldots,f_k \}.$  

We let $\mathcal{C}(\mathbb{R}^n)$ denote the space of compact convex subsets of 
$\mathbb{R}^n$ and endow it with the Hausdorff metric 
$\textup{d}_{\textup{H}}$ given by 
$$\textup{d}_{\textup{H}}(A,\, B):=
\inf\{\epsilon>0 \, \vert \,  A\subset N_{\epsilon}(B)\, \text{and} \, B\subset N_{\epsilon}(A) \},$$ 
where $N_{\epsilon}(X)$ denotes the open $\epsilon$\--neighborhood 
of a subset $X\subset \mathbb{R}^n$.  A compact convex set 
$A\in \mathcal{C}(\mathbb{R}^n)$ with non-empty interior is said to be a \textit{convex body}.  
If $\lambda>0$ and $A$ and $B$ are convex bodies then so are the sets 
$$\lambda \, A:=\{\lambda \, a \, \vert \, a 
\in A\}\,\,\,\,\,\,\,\,\,\,\,\,\,\,\,\,\,\,\, A+B=\{a+b \, \vert \, a\in A\, \text{ and }\, b\in B\}.$$  
Subsets 
$A,\, B \subset \mathbb{R}^n$ are said to be \textit{homothetic} if there 
exists $v\in \mathbb{R}^n$ and $\lambda>0$ such that $\lambda \, A + \{v\} = B$.  
We recall the celebrated 
Brunn\--Minkowski inequality (see \cite{Ga} for a detailed survey):  

\begin{Theorem}[Brunn-Minkowski]\label{B-M}
Let $A,\, B$ be convex bodies in $\mathbb{R}^n$ and $0<\lambda <1.$  
Then $$\vol_{\textup{euc}}^{1/n}((1-\lambda)\, A+\lambda \, B)
\geq (1-\lambda)\, \vol_{\textup{euc}}^{1/n}(A)
+ \lambda \, \vol_{\textup{euc}}^{1/n}(B),$$ with equality if and only if $A$ and $B$ are homothetic.
\end{Theorem}

\vskip 5pt

In the remainder of this section, we let 
$\Delta^n= \bigcap_{i=1}^{F} \{x \in \mathbb{R}^n \, \vert \,  \langle x,\, u_i \rangle \geq \lambda_i\}$ denote an $n$\--dimensional Delzant polytope with $F$ facets and $V$ vertices.  
We enumerate the vertices $v_1,\, v_2,\ldots,v_V$ and facets 
$\mathcal{F}_1,\ldots, \mathcal{F}_F$.  When two vertices $v_j$ and $v_k$ are 
adjacent in $\Delta$, we denote their common edge by $e_{j,k}$.  By the Delzant property, each vertex $v_i$ is the unique intersection point of $n$ facets, 
$$\mathcal{F}_{v_i}^j=\{x\in\mathbb{R}^n \, \vert \, \langle x,\, u_{v_i}\rangle=
\lambda_{v_i}\}\cap\Delta,\,\,\,\, j=1,\ldots,n,$$ where the inward normal vectors 
$u_{v_i}^j$ to the $j^{\textup{th}}$ facet $\mathcal{F}_{v_i}^j$ collectively define a 
$\mathbb{Z}$\--basis of $\mathbb{Z}^n$.  In particular, $\Delta$ is simple and 
there are precisely $n$ edges $e_{v_i}^j$, $j=1,\ldots n$ leaving each 
vertex $v_i$ of $\Delta$.  We shall denote the set of admissible packings 
of $\Delta$ by $\mathcal{AP}(\Delta)$.

\begin{named}{Proposition \ref{convex prop2}}
Let $\Delta \subset \mathbb{R}^n$ be a Delzant polytope. The set 
$\mathcal{AP}(\Delta)$ of admissible packings
of $\Delta$ has the structure of a polytope in $\mathbb{R}^V$, where $V$
is the number of vertices of $\Delta$.  Moreover, $\Delta$ admits finitely many maximal toric ball packings.
\end{named}

\begin{proof}
Each packing $\mathcal{Q} \in \mathcal{AP}(\Delta)$ 
consists of a disjoint union $\bigsqcup_{i=1}^{V} \Sigma(v_i,\, R_{i}(\mathcal{Q}))$ of 
admissible simplices $\Sigma(v_i,\, R_{i}(\mathcal{Q}))$ centered at the 
vertex $v_{i}$ with (possibly zero) radius $R_{i}(\mathcal{Q})$.  Define the 
map   $$R:\mathcal{AP}(\Delta) \rightarrow \Pi_{i=1}^{V} [0,\, r_{v_i}],\,\,\,\, \mathcal{Q} 
\mapsto (R_{1}(\mathcal{Q}),\ldots, R_{V}(\mathcal{Q})).$$

By Lemma \ref{adsi}, the map $R$ is injective so that we can identify $\mathcal{AP}(\Delta)$ with its image in $\Pi_{i=1}^{V} [0,\, r_{v_i}]$.  Note that $R$ is not surjective as the admissible simplices with given radii  $(x_1,\ldots, x_V) \in  \Pi_{i=1}^{d} [0,\, r^{v_i}]$ will not in general be disjoint.  We must argue that the image set of $R$ is precisely the solution set to finitely many linear inequalities. 

As a first step, we give a criterion for admissible simplices to be disjoint in terms of their intersections along the edges of $\Delta$.
To this end, let $\mathcal{F}$ denote a finite family of admissible simplices with pairwise distinct 
centers and fix an admissible simplex $\Sigma$ from the family.  Let $v$ denote the center 
of $\Sigma$.    

\vskip 10pt
{\bf \textup{Disjointness Condition}}: \textit{For $\Sigma$ to be disjoint from the rest of the family $\mathcal{F}$, it is a necessary and sufficient condition that its closure $\overline{\Sigma}$ intersects the closure of the other admissible simplices in the family in at most one point in each of the edges $e_{v}^j$, $j=1,\ldots, n$.}  

\vskip 10pt

To see that the above disjointness condition is necessary, suppose that the closure of another admissible simplex in $\mathcal{F}$, say $\overline{\Sigma'}$, intersects $\overline{\Sigma}$ in more than a single point along some edge $e_{v}^j$. Then $\overline{\Sigma}\cap\overline{\Sigma'}\cap e_v^j$ is a convex subset of $e_v^j$ with at least two points and is therefore a closed subsegment $e\subset e_i^v$ with nonempty interior.  The interior of $e$ is contained in both $\Sigma$ and $\Sigma'$, whence $\Sigma \cap \Sigma' \neq \emptyset$.   Next we argue that the above condition is sufficient.  To see this, let $\Sigma'$ be another admissible simplex in the family $\mathcal{F}$. Denote by $v'\in \Delta$ the center of $\Sigma'$ and let $x_{v'}^j$ denote the unique point in the set 
$\overline{\Sigma'}\cap e_{v'}^j \setminus \Sigma' \cap e_{v'}^j$ for each $j=1, \ldots, n$.  If the above condition holds, then $\{v',\,x_{v'}^1, \ldots , x_{v'}^n\} \subset \Delta \setminus \Sigma$.  But since $\Delta \setminus \Sigma$ is a convex set and $\overline{\Sigma'}= \ConvHull(v',\,x_{v'}^1, \ldots , x_{v'}^n)$, it follows that $\overline{\Sigma'}\subset \Delta \setminus \Sigma,$ concluding the proof of sufficiency.
\vskip 10pt

For distinct vertices $v_i,\, v_j \in \{v_1,\ldots, v_V\}$, define $L_{i,j}$ by
$$
L_{i,j}= \left\{ \begin{array}{rl}
         \length_{\mathbb{Q}}(e_{i,j}) & \mbox{ if } v_i\, \mbox{and}\, v_j\, \mbox{are adjacent}\\
      r_{v_i}+r_{v_j}  & \mbox{ otherwise } \end{array} \right .
$$

It follows from (\ref{rv}) and Lemma \ref{adsi} and the preceding disjointness condition that a point $(x_1,\ldots,x_V) \in \Pi_{i=1}^{V} [0,\, r_{v_i}]$ lies in the image of $R$ if and only if the equations 
$$x_{i}+x_{j} \leq L_{i,j},$$ hold for all $i \neq j \in \{1,\ldots, V\}$.  Therefore,
\begin{eqnarray} \label{ap}
R(\mathcal{AP}(\Delta))= 
\bigcap_{i\neq j\in\{1,\ldots,V\}} \{(x_1,\ldots,x_V)\in \mathbb{R}_{\geq0}^V\,\vert\, x_i+x_j\leq L_{i,j}\}.
\end{eqnarray}

By Lemma \ref{adsi}, the density function in the coordinates of $\mathbb{R}_{\geq 0}^V$ is given by $$\Omega((x_1,\ldots, x_V))=(\Sigma_{i=1}^V x_{i}^{n})/(n!\, \vol_{\textup{euc}}(\Delta)).$$  As $n>1$, $\Omega$ is a strictly convex function on $\mathcal{AP}(\Delta)$ so that its maximum value can only be obtained at its vertices,  establishing the last part of the proposition.
\end{proof}

In view of the injectivity of $R$, we will henceforth identify $\mathcal{AP}(\Delta)$ with the polytope $R(\mathcal{AP}(\Delta)).$

Before proving Proposition \ref{main prop}, we need to set up some additional notation and prove a lemma.  A natural way to perturb the Delzant polytope $\Delta$ is by perturbing its defining linear equations.  Define the map 
\begin{eqnarray} \label{rfs}
\mathbb{R}^F \rightarrow \mathcal{C}(\mathbb{R}^n),\,\,\,\,\,\,\, s=(s^1,\ldots,s^F) \mapsto \Delta_{s}
\end{eqnarray}
by $$\Delta_s:= \bigcap_{i=1}^{F} \{x \in \mathbb{R}^n \, \vert \, 
\langle x,\, u_i \rangle \geq \lambda_i +s^i\}.$$  This map is continuous and has image in the set of polytopes with not more than $F$ facets (we regard the empty set as a polytope). For small 
enough $r>0$, the polytopes $\Delta_s$ with $s\in \mathbb{B}(0,\, r)\subset \mathbb{R}^F$ still have $F$ facets, are Delzant, and all determine the same fan. Let $B_{\Delta}$ denote the largest open ball centered at the origin with these three properties and let $\mathcal{D}^n$ denote the set of equivalence classes of Delzant polytopes.  The map $B_{\Delta} \rightarrow \mathcal{D}^n$ induced by (\ref{rfs}) induces an equivalence relation on $B_{\Delta}$ by declaring points in the fibers to be equivalent.  We denote the set of equivalence classes by $\widehat{B}_{\Delta}$ and endow it with the quotient topology.  As $\textup{GL}(n,\,\mathbb{Z})$ is discrete, the dimension of a fixed equivalence class of Delzant polytopes is $n+1$, the dimension of the group of homotheties of $\mathbb{R}^n$.  It follows that the dimension of $\widehat{B}_{\Delta}$ satisfies $\dim(\widehat{B}_{\Delta})\geq F - (n+1).$ 

Define the function  $$\Omega:B_{\Delta} \rightarrow (0,1]$$ by $\Omega(s):=\Omega(\Delta_s)$, and  for $s_1,\, s_2 \in B_{\Delta}$, define the function 
$$
\Omega_{s_1,\,s_2}:[0,\, 1] \rightarrow (0,\,1]
$$ 
by $\Omega_{s_1,\, s_2}(t):=\Omega((1-t)\, s_1+t\, s_2)=\Omega(\Delta_{(1-t)\, s_1
+t\, s_2})$.

\begin{lemma}\label{regularity}
The function $\Omega:B_{\Delta} \rightarrow (0,\,1]$ is continuous.  For distinct $s_1,\, s_2 \in B_{\Delta}$, there exists a suitably small $\epsilon(s_1,\, s_2)>0$ such that the restriction $\Omega_{s_1,\,s_2}^{1/n}|_{[0,\, \epsilon]}:[0,\, \epsilon] \rightarrow (0,\,1]$ is convex.  Furthermore, if $\Delta_{s_1}$ and $\Delta_{s_2}$ are not homothetic, then $\Omega_{s_1,\, s_2}^{1/n}|_{[0,\, \epsilon]}$ is strictly convex. 
\end{lemma}

\begin{proof}
Define $||\cdot||_n:\mathbb{R}_{\geq 0}^V \rightarrow \mathbb{R}$ by $||(x_1,\ldots,x_V)||_n=(\sum_{i=1}^V x_i^n)^{1/n}$.  It is well known that $||\cdot||_n$ is a strictly convex function.  For $s\in B_{\Delta}$, we have that 
$$\Omega^{1/n}(s)=
\frac{\max ||\cdot||_n|_{\mathcal{AP}(\Delta_s)}}{(n!)^{1/n}\vol_{\textup{euc}}^{1/n}(\Delta_s)}.$$

Each of the vertices $v_1,\ldots,v_V$ of $\Delta$ are defined as the unique solution to the linear system:
$$\langle u_{v_i}^j,\, v_i \rangle=\lambda_{v_i}^j\,\,\,\,\,\,\, j=1,\ldots n.$$  
Similarly, $\Delta_s$ has $V$ vertices $v_1(s),\ldots,v_V(s)$, each defined by the linear system:
$$\langle u_{v_i}^j,\, v_i(s) \rangle=\lambda_{v_i}^j+s_{v_i}^j\,\,\,\,\,\,\, j=1,\ldots n.$$

Hence, the vertices define affine maps $v_i:B_{\Delta} \rightarrow \mathbb{R}^n$.  It follows that the edges of $\Delta_s$ and their rational lengths also vary linearly with $s\in B_{\Delta}$ and that, as we
have already remarked, $s \mapsto \Delta_s$ defines a continuous map $B_{\Delta}\rightarrow \mathcal{C}(\mathbb{R}^n)$. 

By the proof of Proposition \ref{convex prop2}, c.f. expression (\ref{ap}), 
$$\mathcal{AP}(\Delta_s)= \bigcap_{i\neq j\in\{1,\ldots,V\}} \{(x_1,\ldots,x_V)\in \mathbb{R}_{\geq0}^V\,\vert\, x_i+x_j\leq L_{i,j}(s)\},$$ 
where 
$$
L_{i,j}(s) = \left\{ \begin{array}{rl}
         \length_{\mathbb{Q}}(e_{i,j}(s)) & \mbox{ if } v_i(s)\, \mbox{and}\, v_j(s)\, \mbox{are adjacent};\\
      r_{v_i}(s)+r_{v_j}(s) & \mbox{ otherwise, } \end{array} \right .
$$            
and 
$$r_{v_i}(s)=\min\{\length_{\mathbb{Q}}(e_{v_i}^1(s)),\ldots, \length_{\mathbb{Q}}(e_{v_i}^n(s)) \}.$$  Hence, the defining equations of $\mathcal{AP}(\Delta_s)$ vary continuously with $s\in B_{\Delta}$ so that $s \mapsto \mathcal{AP}(\Delta_s)$ defines a continuous map $B_{\Delta} \rightarrow \mathcal{C}(\mathbb{R}_{\geq 0}^V)$.  Since $\vol^{1/n}$ and $\max ||\cdot||_n$ define continuous maps on $\mathcal{C}(\mathbb{R}^n)$ and $\mathcal{C}(\mathbb{R}_{\geq 0}^V)$, it follows that $\Omega$ is continuous.

Now fix $s_1,\, s_2 \in B_{\Delta}$ and $t \in (0,\,1)$.  We first claim that 
$$\Delta_{(1-t)\, s_1+t\, s_2}=(1-t)\, \Delta_{s_1}+t\, \Delta_{s_2}.$$  
Note that 
$$\Delta_{(1-t)\, s_1+t\, s_2}=\ConvHull(v_1((1-t)\, s_1+t\, s_2),\ldots,v_V((1-t)\, s_1+t\, s_2)).$$  
Since $v_{i}((1-t)\, s_1+t\, s_2)=(1-t)\, v_i(s_1)+t\, v_i(s_2),$ 
$$
\{v_1((1-t)\, s_1+t\, s_2),\ldots,v_V((1-t)\, s_1+t\, s_2)\}\subset (1-t)\, \Delta_{s_1}
+t\, \Delta_{s_2},
$$ whence 
$$
\Delta_{(1-t)\, s_1+t\, s_2}\subset (1-t)\, \Delta_{s_1}+t\, \Delta_{s_2}.
$$  
Now consider the $(n+1)$-dimensional polytope 
$$\Delta(s_1,\, s_2):=\bigcap_{i=1}^{F}\{(x,\, t) \in \mathbb{R}^n\times[0,\,1] \, \vert \, \langle x,\, u_i \rangle \geq \lambda_i +(1-t)\, s_1+t\, s_2\}.
$$
Let $H_t=\Delta(s_1,\, s_2) \cap \{x_{n+1}=t\}$ and note that $H_t$ is naturally identified with $\Delta_{(1-t)\, s_1+t\, s_2}$.  Now, if $(x,\, 0)\in H_0$, $(y,\, 1) \in H_1$, and 
$t\in (0,\,1)$, then $(1-t)\, (x,\, 0)+t\, (y,\, 1) \in H_t$, concluding the proof of the claim.

By Theorem \ref{B-M}, the map $[0,\, 1] \rightarrow \mathbb{R}$ given by 
$t \mapsto \vol_{\textup{euc}}^{1/n}(\Delta_{(1-t)\, s_1+t\, s_2})$ is concave and strictly concave if and only if  $\Delta_{s_1}$ and $\Delta_{s_2}$ are not homothetic.  To conclude the proof of the Lemma, we must argue that there is a suitably small $\epsilon=\epsilon(s_1,\, s_2)>0$ such that the map $[0,\,1]\rightarrow \mathbb{R}$ given by $t \mapsto \max ||\cdot||_n|_{\mathcal{AP}(\Delta_{(1-t)\, s_1+t\, s_2})}$ is convex when restricted to the interval $[0,\, \epsilon)$.  Let $p_1,\, p_2,\ldots,p_k \in \mathbb{R}^V$ denote the vertices of $\mathcal{AP}(\Delta_{s_1})$.  It follows from the description above that for $t$ sufficiently small $\mathcal{AP}(\Delta_{(1-t)\, s_1+t\, s_2})$ also has $k$ vertices $v_1(t),\ldots,v_k(t)\in \mathbb{R}^V$ and that the maps $t\mapsto v_i(t)$ define (possibly constant) line segments in $\mathbb{R}^V$.  By convexity of $||\cdot||_n$, $$\max ||\cdot||_n|_{\mathcal{AP}(\Delta_{(1-t)\, s_1
+t\, s_2})}=\max \{||v_1(t)||_n,\ldots,||v_k(t)||_n\},$$ and the result follows.
\end{proof}

We are now ready to prove:

\begin{named}{Proposition \ref{main prop}}

Let $\mathcal{D}^n$ denote the set of equivalence classes of $n$\--dimensional Delzant polytopes and $$\Omega:\mathcal{D}^n\rightarrow (0,\,1]$$ be the maximal density function, $n\geq2$.  
Suppose that $\Delta$ is a Delzant polytope having at least $n+3$ facets and let $\Omega(\Delta):=\delta \in (0,\,1)$.
Then for any $\epsilon>0$,  
there exists a constant $c>0$ and a family $\mathcal{F}$ of Delzant polytopes satisfying 

\begin{itemize}

\item 
the polytopes from $\mathcal{F}$ determine a common fan,

\item 
$|\vol_{\textup{euc}}(\Delta')-\vol_{\textup{euc}}(\Delta)| < \epsilon$ for all $\Delta'  \in \mathcal{F}$,

\item $\Omega^{-1}(\{x\}) \cap \mathcal{F}$ is uncountable for all
$x \in (\delta-c,\, \delta)$ or for all $x \in (\delta, \, \delta+c)$. 

\end{itemize}

\end{named}

\begin{proof}
Let $\Delta$ be a Delzant polytope with with at least $n+3$ facets and with $\Omega(\Delta)=\delta$ and let $\epsilon>0$.  By continuity of the volume function with respect to the Hausdorff metric, there is a suitably small connected neighborhood $N\subset B_{\Delta}$ of the origin for which $|\vol_{\textup{euc}}(\Delta_s)-\vol_{\textup{euc}}(\Delta)|<\epsilon$ for each $s \in N$.  We define the desired family $\mathcal{F}$ by $$\mathcal{F}=\{\Delta_s \, \vert \, \, s \in N\}$$ and remark that by construction all $\Delta' \in \mathcal{F}$ determine the same fan.  Therefore, it remains to establish the third item of the proposition.

As $\textup{dim}(N) \ge n+3$ and since the space of homotheties of $\mathbb{R}^n$ has dimension $n+1$, there exists $s \in N \setminus \{0\}$ such that $\Delta_s$ is not
homothetic to $\Delta$. By Lemma \ref{regularity}, there exists 
$\epsilon'>0$ such that $\Omega^{1/n}_{0,\,s}\colon [0, \, \epsilon') \to [0,\,1)$
is strictly convex and therefore
$\Omega \colon N \to (0,\,1]$ is not the constant map.  Since
$\Omega$ is continuous and $N$ is connected, there exists $c>0$, 
such that $(\delta, \, \delta+c) \subset \Omega(N)$ or
$(\delta-c, \, \delta) \subset \Omega(N)$.

Suppose that $(\delta-c,\,\delta) \subset \Omega(N)$.  Note that $\Omega:B_{\Delta} \rightarrow (0,\,1]$ descends to a continuous map $\widehat{\Omega}:\widehat{B}_{\Delta}\rightarrow (0,\,1]$.  Suppose that for some $x \in (\delta-c,\, \delta)$, $\widehat{\Omega}^{-1}(\{x\})$ is countable.  As $\dim(\widehat{B}_{\Delta})\geq 2$, $\widehat{B}_{\Delta}\setminus \widehat{\Omega}^{-1}(\{x\})$ is connected.  Hence,  $\widehat{\Omega}(\widehat{B}_{\Delta}\setminus \widehat{\Omega}^{-1}(\{x\}))$ is connected, a contradiction.  
In case that $(\delta, \, \delta+r) \subset \Omega(N)$ an analogous argument
concludes the proof.
\end{proof}

\begin{remark} \label{rmk}
\normalfont
The proof of Proposition \ref{main prop} above also establishes the following:  suppose that $\Delta$ is an $n$-dimensional Delzant polytope with at least $n+3$ facets and with $\Omega(\Delta):=\delta \in (0,\,1)$.  Also suppose that $\Omega(B_{\Delta})$ contains an open neighborhood $(\delta-c,\, \delta+c)$ of $\delta$.  Then for each $x \in (\delta-c,\, \delta+c)$, $\Omega^{-1}(\{x\}) \subset \mathcal{D}^n$ is uncountable.
\end{remark}

\vskip 5pt

\begin{named}{Proposition \ref{contrast prop}}
Let $\mathcal{D}^n$ denote the set of equivalence classes of $n$\--dimensional Delzant polytopes and $$\Omega:\mathcal{D}^n\rightarrow (0,\,1]$$ be the maximal density function, $n\geq 2$.  Then $\Omega^{-1}(\{x\})$ is uncountable for all $x \in (0,\, 1)$. 
\end{named}

\begin{proof}
By Remark \ref{rmk}, it suffices to show that there exists an $n$-dimensional Delzant polytope $\Delta$ with at least $n+3$ facets for which $\Omega(B_{\Delta})=(0,\,1)$.  Consider the polytope $\Delta(\epsilon_1,\,\epsilon_2)$ obtained by removing from the standard $n$\--dimensional simplex an admissible simplex of radius $\epsilon_i$ at the vertex $e_i$ for $i=1,\,2$.  For compatible choices of $\epsilon_1$ and $\epsilon_2$, we obtain a Delzant polytope with $n+3$ facets.  Fix $\epsilon_1^0$ and $\epsilon_2^0$ both very close to zero and  let $E\subset [0,\,1]^2$ be the set of pairs $(x,y)$ for which $\Delta(x,\, y)=\Delta_{s}$ for some $s\in B_{\Delta(\epsilon_1^0,\, \epsilon_2^0)}$.  By definition, $\Omega(\{\Delta(x,\,y)\, \vert \, (x,\,y) \in E\}) \subset \Omega(B_{\Delta(\epsilon_1^0,\,\epsilon_2^0)})$, a connected subset of the interval $(0,\,1)$.  We conclude by showing $\Omega(\{\Delta(x,\,y)\vert (x,\,y) \in E\})$ contains deleted open neighborhoods of $0$ and 
$1$ in $[0,\,1]$.  To obtain a deleted neighborhood of $1$ we remark that as $(x,\,y) \rightarrow (0,\,0)$, $(x,\,y) \in E$ and $\Omega(\Delta(x,\,y)) \rightarrow 1$.  Similarly, to obtain a deleted neighborhood of $0$ we remark that as $y \rightarrow 0$, the pairs of the form $(1-2y,\,y)$ are in $E$ and $\Omega(1-2y,\,y) \rightarrow 0$.  
\end{proof}

We thank Professor Novik for bringing the example above to our attention.  We conclude with the following:

\vskip 5pt

\textbf{Question.} Let $[(M^{2n},\,\sigma, \, \psi)]$ be an equivalence class of symplectic\--toric manifolds.  Is there a formula for the number of different maximal toric packings of $M$ in terms of its  equivariant symplectic invariants?

\noindent
Alvaro Pelayo\\
Department of Mathematics, University of Michigan\\
2074 East Hall, 530 Church Street, Ann Arbor, MI 48109--1043, USA\\
e\--mail: apelayo@umich.edu

\bigskip

\noindent
Benjamin Schmidt\\
Department of Mathematics, University of Chicago\\
5734 South University Avenue,
Chicago, Illinois 60637\\
e\--mail: schmidt@math.uchicago.edu

\end{document}